\newcommand{\E}[0]{{\sf E}}
\newcommand{\one}[0]{{\bf 1}}
\newcommand{\uzero}[0]{\underline{0}}
\newcommand{\uone}[0]{\underline{1}}

\newcommand{\nn}[0]{{\bf n}}
\newcommand{\mm}[0]{{\bf m}}

\documentclass[letterpaper,11pt]{article}

\usepackage{amsthm}
\usepackage{amsmath}

\newtheorem{thm}{Theorem}[section]

\newtheorem{lemma}[thm]{Lemma}

\newtheorem{question}[thm]{Question}

\newcommand{\beq}[1]{\begin{equation}\label{#1}}
\newcommand{\enq}[0]{\end{equation}}

\newcommand{\bn}[0]{\bigskip\noindent}
\newcommand{\mn}[0]{\medskip\noindent}
\newcommand{\nin}[0]{\noindent}

\newcommand{\sub}[0]{\subseteq}
\newcommand{\sm}[0]{\setminus}
\renewcommand{\dots}[0]{,\ldots,}

\newcommand{\f}[0]{{\cal F}}

\newcommand{\ra}[0]{\rightarrow}
\newcommand{\Ra}[0]{\Rightarrow}

\newcommand{\ff}[0]{{\bf f}}

\newcommand{\TTT}[0]{R}
\newcommand{\SSS}[0]{T}

\newcommand{\0}[0]{\emptyset}

\renewcommand{\qed}[0]{\begin{flushright} \rule{2mm}{3mm} \end{flushright}}

\newcommand{\C}[2]{{{#1}\choose{{#2}}}}
\newcommand{\Cc}[0]{\tbinom}
\newcommand{\ga}[0]{\alpha }
\newcommand{\gb}[0]{\beta }

\newcommand{\gd}[0]{\delta }

\newcommand{\go}[0]{\omega}
\newcommand{\gO}[0]{\Omega}

\newcommand{\gz}[0]{\zeta}
\newcommand{\eps}[0]{\varepsilon }

\newcommand{\vp}[0]{\varphi}

\newcommand{\comments}[1]{}
\usepackage[usenames]{color}
\begin{document}
\renewcommand{\thefootnote}{\fnsymbol{footnote}}
\footnotetext{AMS 2010 subject classification:  }
\footnotetext{Key words and phrases:  }

\title{Functions without influential coalitions}
\author{J. Kahn\footnotemark $~$ and G. Kalai\footnotemark
}
\date{}
\footnotetext{* Supported by NSF grant DMS0701175  and BSF grant 2006066.}
\footnotetext{ $\dag$ Supported by BSF grant 2006066. }

\date{}

\maketitle

\begin{abstract}

We give counterexamples to a conjecture of Benny Chor
and another of the second author, both from the late 80s,
by exhibiting functions for which the influences of
large coalitions are unexpectedly small relative to the
expectations of the functions.

\end{abstract}

\section {Introduction}\label{Intro}

\mn

For a set $T$ we use $\Omega(T)$ for the discrete cube
$\{0,1\}^T$ and $\mu_T$ for the uniform probability measure
on $\Omega (T)$.
In this paper
$f$ will always be a Boolean function on $\gO([n])$
(that is, $f:\gO([n]) \to \{0,1\}$, where, as usual,
$[n]=\{1\dots n\}$),
and we will write $\mu$ for $\mu_{[n]}$.
We reserve $x,y$ for elements of $\Omega([n])$
and set $|x| = \sum x_i$.

Following Ben-Or and Linial \cite {bl} we define,
for a given $f$ and
$S \subset [n]$, the {\em influence of $S$ toward one}
to be
\beq{I+}
I^+_S(f) = \mu_{[n]\sm S}( \{u \in \Omega ([n]\backslash S):
\exists v \in \Omega(S), f(u,v)=1\})- \mu (f).
\enq
Similarly, the influence of $S$ toward zero is
\beq{I-}
I^-_S(f) = \mu_{[n]\sm S}( \{u \in \Omega ([n]\backslash S):
\exists v \in \Omega(S), f(u,v)=0\})- (1-\mu (f))
\enq
and the (total) influence of $S$ is
$$I_S(f)=I_S^+(f)+I_S^-(f).$$

Suppose $\mu(f)=1/2$.
It then follows from a theorem of Kahn, Kalai and Linial
\cite{kkl}
(Theorem \ref{t:kkl} below, henceforth ``KKL")
that for every $a>0$ there is an $S \subset [n]$ of size
$an$ with
$I^+_S(f)\ge 1/2-n^{-c}$, where $c>0$ depends on $a$.
(See Theorem \ref{t:ra1}.)
Benny Chor conjectured in 1989 that one can in fact
achieve
$I^+_S(f)\ge 1/2-c^n$ (where, again, $c<1$ depends on $a$).
The conjecture has been ``in the air" since that time,
though as far as we know
it has appeared in print only in
\cite{k:nati,Nati}.

In this note we disprove
Chor's conjecture and another, similar conjecture
from the same period.

For our purposes the subtracted terms in \eqref{I+} and
\eqref{I-} are mostly a distraction, and it sometimes seems
clearer to speak of $J^+_S(f):=I^+_S(f)+\mu(f)$
and $J^-_S(f):=I^-_S(f)+(1-\mu(f))$.
Thus, for example, $J^+_S(f)$ is the probability that
a uniform setting of the variables in $[n]\sm S$
doesn't force $f=0$, and Chor's conjecture predicts
an $S$ with $J^+_S(f)\ge 1-c^n$.
The following statement shows that this need not be the case.

\begin{thm}\label{prop1}
For any fixed $\ga,\gd \in (0,1)$
there are a $C$ and an $f$ with
$\mu(f)=\ga$ and $J_S^+(f)<1-n^{-C}$
for every $S\sub [n]$ of size $(1/2-\gd)n$.
\end{thm}

\nin
We should note that one cannot expect to go much
beyond $|S|=(1/2-\gd)n$; for example if $\mu(f) =1/2$,
then it follows from the ``Sauer-Shelah
Theorem" (Theorem \ref{t:sas}) that there is
an $S$ of size $n/2$ with
$J_S^+(f) = 1$.

Another consequence of KKL
(see Theorem \ref{t:ra2} below)
is that there is a $\beta>0$ such that
for any
$f$ with $\mu(f)> n^{-\beta}$ there is an
$S$ of size (say) $0.1 n$ with influence $1-o(1)$.
A conjecture of the second author, again from
the late 80s, asserts
that the same conclusion holds even assuming only
$\mu (f)> (1-\eps)^n$
for sufficiently small $\eps$.
This conjecture turns out to be false as well:

\begin{thm}\label{prop2}
For any fixed $\eps, \gd>0$
there are an $\ga >0$ and
Boolean functions $f$ such that
$\mu(f)> (1-\eps)^n$ and no set of size $(1/2-\gd)n$
has influence to 1 more than $\exp[-n^\ga]$.
\end{thm}


\mn
(This can be strengthened a bit to require
$\mu(f) > \exp[-n^{1-c}]$ for some fixed $c=c_\gd>0$.)

It would be interesting to see
how close one can come to the positive result
(i.e. Theorem \ref{t:ra2})
mentioned before Theorem \ref{prop2}.
For example, could it be that there is some fixed $\gb$
for which one can find $f$'s with
$\mu(f)\approx n^{-\gb}$ for which no $S$
of size $0.1 n$ has influence $\gO(1)$?
We will discuss this question further in the next section.

Each of our examples
is of the form $f= \wedge_{i=1}^m C_i$, where the
$C_i$'s are random $\vee$'s of $k$ literals
using $k$ distinct variables (henceforth ``$k$-clauses").
These $f$'s, which may be thought of as
variants of the ``tribes"
construction of Ben-Or and Linial
(see below),
were inspired by a paper of
Ajtai and Linial \cite {AL} and share with it the following
curious feature.
It's easy to see that any $f$ can be converted to a
{\em monotone}
(i.e. increasing) $f'$ with $\mu(f')=\mu(f)$ and each
influence ($I_S^+$ and so on) for $f'$ no larger than
the corresponding influence for $f$;
thus it's natural to look for $f$'s with small influences
among the increasing functions.
But the present random examples, like those of \cite{AL},
do not do this, and it's not easy to see what
one gets by monotonizing them.

\section {Background and perspective}
\label{s:bac}


\subsubsection *{Influence}


We write $I_\ell (f)$ for $I_{\{\ell\}}(f)$.  A form of
the classic edge isoperimetric inequality for Boolean
functions is

\begin{thm}\label{Tiso}
For any (Boolean function) f with $\mu(f)=t$,
\begin {equation}
\label {e:harper}
I(f) := \sum_{k=1}^n  I_\ell (f) \ge  2t \log_2 (1/t).
\end {equation}
\end {thm}

\mn
(This convenient version is easily derived from the precise
statement, due to Hart \cite{Hart}; see also \cite[Sec. 7]{KK07}
for a simple inductive proof.)

\mn
While (\ref{e:harper}) is exact or close to exact
(depending on $t$), it typically gives
only a weak lower bound on the maximum of the $I_\ell (f)$'s,
namely
\begin {equation}
\label {e:harper-max}
\max_\ell I_\ell (f) \ge  2t \log_2 (1/t)/n.
\end {equation}
For $t$ not too close to 0 or 1, the
following
statement from \cite{kkl} gives
better information.

\begin{thm} [KKL]
\label {t:kkl}
There is a fixed $c>0$ such that
for any
$f$ with $\mu(f)=t$,
there is a $k\in [n]$ with
\begin {equation}
\label {e:kkl}
I_\ell (f) \ge c t (1-t)\log n/n.
\end {equation}
\end {thm}
\nin
Recall that $J_\ell (f)=\mu(f)+I_\ell (f)$. Repeated application of Theorem \ref{t:kkl} gives
the following two corollaries.

\begin {thm}
\label {t:ra1}
For all $a, t\in (0,1)$ there is a c such that
for any
$f$ with $\mu(f)=t$ there
is an $S\sub [n]$ with $|S|\leq an$ and
$$J^+_S(f) \ge 1 - n^{-c}$$ (that is,
$I^+_S(f) \ge (1-t) - n^{-c}$).

\end {thm}
\nin
Similarly (either by the same argument or by applying
Theorem \ref{t:ra1} to the function $1-f(x)$)
there is a small $S'$
with
$J^-_{S'}(f) \ge 1 - n^{-c}$ (i.e. $I^-_{S'}(f) \ge t - n^{-c}$),
and combining these observations we find that there is
in fact a small $S''$ (e.g. $S \cup S'$) with
$I_{S''}(f) \ge 1-n^{-c}$.

\begin {thm}
\label {t:ra2}
For every $\delta, \epsilon >0$,
there is an $\ga>0$ such that for large enough
$n$ and any
$f$ with $\mu(f)\geq  n^{-\ga}$,
there is an $S\sub [n]$ with $|S|=\delta n$ and
 $$J^+_S(f) \ge 1-\epsilon.$$
\end {thm}


\mn
The conjecture of
Chor stated in Section \ref{Intro}
asserts that the $n^{-c}$ in
Theorem \ref {t:ra1} can be replaced
by something exponential in $n$,
and the conjecture stated before
Theorem \ref{prop2} proposes a similar weakening
of the $n^{-\ga}$ lower bound on $\mu(f)$ in
Theorem \ref{t:ra2}.
As already noted, we
will show below that these conjectures are incorrect.

\subsubsection *{Tribes}
The original ``tribes" examples of Ben-Or and Linial \cite{bl}
are Boolean functions of the form
$f= \vee_{i=1}^m C_i$, where the
``tribes" $C_i$ are $\wedge$'s of $k$ (distinct) variables
and
each variable belongs to exactly one tribe.
The dual of such an $f$ (so ``dual tribes") is
$g= \wedge_{i=1}^m D_i$, where
$D_i$ is the $\vee$ of the variables in $C_i$
(so again, each variable belongs to exactly one $D_i$).

When $k=\log n - \log \log n - \log \ln (1/t)$,
we have $1-\mu(f) = \mu(g) \approx t$
(where $\log = \log_2$ and $f,g$ are as above).
For fixed $t\in (0,1)$ both constructions
show that Theorem \ref{t:kkl}
is sharp (up to the value of $c$).

On the other hand,
when $t=O(n^{-c})$ for a fixed
$c>0$, $f$ shows that (\ref{e:harper-max})
is tight up to a multiplicative constant, depending on $c$;
for example,
$k=2\log n-\log \log n$ gives
$\mu(f) \approx 1/(2n)$ and
$I_\ell(f)\approx 2\log n/n^2=\Theta(\mu(f)\log(1/\mu(f))/n)$
for each $\ell$.
(In contrast, for $\mu(g)\approx 1/n$, we should
take $k= \log n-2\log \log n-1$, in which case
$I_\ell (g)=\Theta (\log^2n/n^2)$ and \eqref{e:harper-max}
is off by a log.)

For $f$ (again, as above) with $\mu(f) \in (\gO(1),1-\gO(1))$,
there are sets of size $\log n$ with large influence towards 1,
while no set of size $o(n/\log n)$ has
influence $\gO(1)$ towards 0.  (The corresponding statement with
the roles of 0 and 1 reversed holds for $g$.)
The Ajtai-Linial construction
mentioned in the introduction shows that there are Boolean functions
$h$ with $\mu(h) \approx 1/2$ and $I_S(h) < o(1)$ for
every $S$ of size $o(n/\log^2n)$.

\subsubsection *{Trace}

We now briefly consider influences from a different point of view.
For a set $X$ let
$2^X = \{S: S \subset X\}$,
${{X}\choose{k}}=\{S \subset X:|S|=k\}$ and
${{X}\choose {< k}}=\{S \subset X:|S| < k\}$.
For $\f\subset 2^X$
and $Y \subset X$,
the {\it trace} of $\cal F$ on $Y$ is
$${\cal F}_{|Y} = \{S \cap Y:
S \in {\cal F}\}.$$

Let $X = [n]$.
The following ``Sauer-Shelah Theorem"
determines, for every $n$ and $m$, the minimum $T$
such that for each $\f\sub 2^X$ of size $T$ there is some
$Y \in \C{X}{m}$ on which the trace of $\f$ is {\em complete,}
meaning ${\cal F}_{|Y}=2^Y$. Such a $Y$
is said to be {\it shattered} by ${\cal F}$.

\begin {thm}[The Sauer-Shelah Theorem]
\label {t:sas}
If ${\cal F} \subset 2^{X}$ and
$|\f| > {{n} \choose {< r}} $,
then $\f$ shatters some $Y\in \C{X}{r}$.
\end {thm}

\nin
That this is sharp is shown by
${\cal F} = {{X} \choose {< r}}$, the {\it Hamming ball}
of radius $r-1$ about
$\0$ with respect to the usual Hamming metric on
$2^{X}\equiv \gO(X)$.

Theorem \ref {t:sas} was proved
around the same time by Sauer \cite {Sa},
Shelah and Perles \cite{Sh}, and Vapnik and Chervonenkis \cite {VC}.
It has many connections, applications and extensions
in combinatorics, probability theory, model theory, analysis, statistics
and other areas.

The connection between traces and influences is as follows.
Let $f$ be a Boolean function on $\gO([n])$
and
${\cal F} =f^{-1}(1)$.
It is easy to see
(identifying $\gO([n])$ and $2^{[n]}$ as usual)
that for
$S\sub [n]$ and $T = [n] \backslash S$,
$$J^+_S(f)=2^{-|T|}|{\cal F}_{|T}| .$$
Thus, in the language of traces, we are interested in
the effect of relaxing ``$\f$ shatters $Y$"
to require only that
${\cal F}_{|Y}$ contain a large fraction of $2^Y$.

The following arrow notation
(e.g. \cite {Bo,Fr}) is convenient.  Write
$(N,n) \to (M,r)$ if every $\f\sub 2^{[n]}$
of size N has a trace of size at least $M$ on some
$S \in \C{[n]}{r}$; for example the Sauer-Shelah Theorem
says $(\Cc{n}{< r}+1,n) \to (2^r,r).$


\medskip
One might hope that Hamming balls would again give the
best examples in our relaxed setting, which would say, for example,
that for $m \le n$,
\begin {equation}
\label{e:tentc1}
(\Cc{n}{< r}+1,n) \to (\Cc{m}{< r}+1,m).
\end {equation}
But \eqref{e:tentc1}, which was first considered by
Bollob\'as and Radcliffe \cite{BR}
and would have implied both of the conjectures
disproved here,
was shown in \cite{BR} to be false
for fixed $k$ and (large) $m=n/2$.
(For $r=n/2$
and $m=n-1$, it fails for the original tribes example
discussed above.)

A consequence of (\ref{e:tentc1})
is that for fixed $\delta, \epsilon>0$ and large $r$,
$$(\Cc{n}{< (1+\epsilon )r/2},n) \to ((1-\delta)2^r ,r),$$
which
would imply our second conjecture from the introduction.
Here a counterexample with $n \gg r$ was given
by Kalai and Shelah \cite{KS}, but this seems not very
relevant to present concerns, for which the
regime of interest has $n$ a little smaller than $2r$.

\subsubsection* {A problem}

\begin{question}\label{prop3}

For fixed $\alpha,\delta >0$, what is
the largest $t\in (0,1/2)$ for which one can find
Boolean functions
$f$ with
$\mu(f)=t$ and
$I^+_S(f)<\ga$ for every $S\sub [n]$ of size
$(1/2-\gd)n$?
\end{question}

\nin
As far as we know $t> n^{-\gb}$ (with $\gb$ depending on $\ga,\gd$)
is possible.



\section{Examples}

In each construction we consider, for suitable $k$ and $m$,
$f= \wedge_{i=1}^m C_i$, where the
$C_i$'s are random $\vee$'s of $k$ literals
using $k$ distinct variables (henceforth ``$k$-clauses")
and show that $f$ is likely to have the desired properties.
We use $g_i$ for the specification associated with $C_i$,
and write $C_i\sim x$ if some entry of $x$ agrees with $g_i$.
We say $C_i$ {\em misses} $S\sub [n]$ if the indices of all
variables in $C_i$ lie in $[n]\sm S$.

Let $s= (1/2-\gd)n$.  We will always
use $S$ for an $s$-subset of $[n]$ and
(for such an $S$) set
$m_S=|\{i:\mbox{$C_i$ misses $S$}\}|.$
(Following common practice we omit irrelevant
floor and ceiling symbols,
pretending all large numbers are integers.
As in the case of $k,m$ and $s$, parameters not declared
to be constants are assumed to be functions of $n$.)
We use $\log$ for $\log_2$.

Both constructions will make use of the next two
observations, with Theorem \ref{prop1}
following immediately from these and
Theorem \ref{prop2} requiring a little more work.
\begin{lemma}\label{mSlemma}
If $k=o(\sqrt{n})$ and $(1/2+\gd)^km=\go(n)$ then
w.h.p.
\beq{mS}
m_S \sim (1/2+\gd)^km ~~\forall S\in \Cc{[n]}{s}.
\enq
\end{lemma}
\nin
(where, as usual, $a_n\sim b_n$ means $a_n/b_n\ra 1 $
and
{\em with high probability} (w.h.p.) means
with probability tending to 1, both as $n\ra\infty$).

\mn
{\em Proof.}
For a given $S$, $m_S$ has the binomial distribution
$B(m,p)$, with $p = \C{n-s}{k}/\C{n}{k}\sim(1/2+\gd)^k$
(using $k=o(\sqrt{n})$ for the ``$\sim $").
Thus
$\E m_S=mp$ and,
by ``Chernoff's Inequality" (e.g.
\cite[Theorem 2.1]{JLR}),
$$
\Pr(m_S\not\in ((1-\gz)mp,(1+\gz)mp) < \exp [-\gO(\gz^2mp)],
$$
for $\gz\in (0,1)$.
Applying this with a $\gz$ which is both
$\go(\sqrt{n/(mp)})$ and $o(1)$ gives
$\Pr(m_S\not\sim mp )< 2^{-\go(n)}$,
and the union bound then gives \eqref{mS}.\qed

\medskip
The next lemma is stated to cover
both applications, though nothing so precise
is needed for Theorem \ref{prop1}.

\begin{lemma}\label{muflemma}
If there is a $\xi$ for which
\beq{xi1}
\exp[-\xi^2n] = o((1-2^{-k})^m)
\enq
and
\beq{xi2}
[(1+2\xi)/4]^k=o(1/m),
\enq
then
w.h.p.
\beq{muf}
\mu(f)\sim (1-2^{-k})^m.
\enq
\end{lemma}

\mn
{\em Proof}.
This is a simple second moment method calculation
(similar to what's done in \cite{AL}, though described
differently there).

Recalling that $x,y$ always denote elements of $\{0,1\}^n$,
write $A_x$ for the event $\{f(x)=1\}$ and $\one_x$
for its indicator,
and set $X=\sum \one_x=2^n\mu(f)$.
Then
$\Pr(A_x)= (1-2^{-k})^m$ and
$\E X = (1-2^{-k})^m2^n$; so we just need to show
$
\E X^2\sim \E^2 X
$ 
(equivalently, $\E X^2<(1+o(1))\E^2 X$),
since Chebyshev's Inequality then gives
$\Pr(|X-\E X|> \gz\E X) =o(1)$
for any fixed $\gz>0$.

We have
$$
\E X^2 = \sum_x\sum_y\E \one_x\one_y =\sum_x\Pr(A_x)\sum_y\Pr(A_y|A_x),
$$
so will be done if we show that for a fixed $x$,
$$ 
\sum_y\Pr(A_y|A_x) <(1+o(1)) (1-2^{-k})^m2^n.
$$ 
Since the sum is the same for all $x$, it's enough
to prove this when $x=\uzero$.
Set $Z=\{y:|y|<(1/2-\xi)n\}$ and recall that by Chernoff's Inequality,
$|Z|<\exp[-2\xi^2n]2^n$.
It is thus enough to show that (for $x=\uzero$)
\beq{nicey}
y\not\in Z ~\Ra~\Pr(A_y|A_x)<(1+o(1)) (1-2^{-k})^m,
\enq
since then, using \eqref{xi1}, we have
$$
\sum_y\Pr(A_y|A_x) ~<~ |Z| + \sum_{y\not\in Z}\Pr(A_y|A_x)
~<~ (1+o(1)) (1-2^{-k})^m2^n.
$$

Now since $x=\uzero$, we have
$A_x =\{g_i\neq \uone ~\forall i\}$; so
if, for a given $y\not\in Z$, we set
$\gb =\gb_y=\Pr(C_i\sim y|g_i\neq \uone)$
(a function of $|y|$),
then $\Pr(A_y|A_x) = \gb^m$.
Aiming for a bound on $\gb$, we have
\begin{eqnarray*}
1-2^{-k}&=&\Pr(C_i\sim y) \\
&=&
\Pr(g_i=\uone)\Pr(C_i\sim y|g_i=\uone)
+ \Pr(g_i\neq\uone)\Pr(C_i\sim y|g_i\neq \uone)\\
&=& 2^{-k}\Pr(C_i\sim y|g_i=\uone) +(1-2^{-k})\gb
\end{eqnarray*}
and
$$\Pr(C_i\sim y|g_i=\uone) > 1 - (1-|y|/n)^k
\geq 1-(1/2+\xi)^k=:1-\nu
$$
(using the fact that if
$g_i=\uone$, then $g_i\not\sim y$
iff all indices of
variables in $C_i$ belong to $\{j:y_j=0\}$).
Combining, we have
$$
\gb < (1-2^{-k})^{-1}[1-2^{-k} - 2^{-k}(1-\nu)]
=(1-2^{-k})\left[1+\tfrac{2^{-k}(\nu -2^{-k})}{(1-2^{-k})^2}\right],
$$
which with \eqref{xi2} gives
$\gb^m < (1+o(1))(1-2^{-k})^m$ (which is \eqref{nicey}).\qed

\mn
{\em Proof of Theorem} \ref{prop1}.
Notice that it's enough to prove this with $\mu(f)\sim \ga$
(rather than ``$=\ga$");
for then,
since $f^{-1}(1)\sub g^{-1}(1)$ trivially implies
$J_S^+(f)\leq J_S^+(f)$ for all $S$, we can
choose
$\gb\in (\ga,1)$ and a $g$ with $\mu(g)\sim \gb$
possessing the desired small influences, and
shrink $g^{-1}(1)$ to produce $f$.

Let $k=C\log n$, with $C=C_\gd$ chosen so that
$(1+2\gd)^k=\go(n)$ (e.g. $C=1/\gd$ does this),
and
$m = 2^k \ln (1/\ga)=n^C\ln (1/\ga)$.
Here all we use from Lemma \ref{mSlemma}
(whose hypotheses are satisfied for our choice of $k$ and $m$)
is the fact that
w.h.p. $m_S\neq 0$ for all $S$,
whence each $J_S^+(f)$ is at most
$1-2^{-k}=1-n^{-C}$.
%
On the other hand, by Lemma \ref{muflemma}
(with, for example, $\xi=0.1$), we have
$\mu(f)\sim \ga$ w.h.p.
So w.h.p. $f$ meets our requirements.\qed

\mn
{\em Proof of Theorem} \ref{prop2}.
Here, intending to recycle $n$, $m$ and $f$, we rename
these quantities
$\nn$, $\mm$ and $\ff$.
We may of course assume $\gd$ is fairly small.
Let (for example) $\xi = \gd/3$, fix $\eps$ with
$0<\eps < \xi^2$, and set
$k=(1+\gd)\log \nn$ and $\mm = \eps 2^k \nn$.
These values are easily seen to give the hypotheses of
Lemmas \ref{mSlemma} and \ref{muflemma}.
In particular, we can say that w.h.p.
the supports of the $C_i$'s
are chosen so \eqref{mS} holds
(note this says nothing about the values specified by
the $g_i$'s)
and
\beq{mS'}
\mu(\ff)\sim (1-2^{-k})^\mm\sim e^{-\eps \nn}.
\enq

Set $n=(1/2+\gd) \nn$.
Fix $S$ ($\in \C{[\nn]}{s}$), set $m=m_S$,
and
let $f=f_S$ be the $\wedge$ of the
$m$ $C_i$'s---w.l.o.g. $C_1\dots C_m$---that miss $S$.
Thus $f$ is the $\wedge$ of
$m\sim (1/2+\gd)^k\mm =\eps (1+2\gd)^k\nn$
random $k$-clauses from a universe of $n$ variables.
Theorem \ref{prop2} (with $\ga=\gd$) thus follows from

\mn
{\em Claim A.}
$\Pr(\mu(f) > \exp[-n^{\gd}]) < o(2^{-\nn})$

\mn
(since then w.h.p. we have $\mu(f_S)\leq \exp[-n^{\gd}]$ for every $S$).

\mn
{\em Remarks.}
The actual bound in Claim A will be
$\exp[-\gO(m)]$,
so much smaller than $2^{-\nn}$.
Note that here it doesn't matter whether we take $\mu$ to be
our original measure (i.e. $\mu$ uniform on $\{0,1\}^\nn$)
or uniform measure on $Q:=\{0,1\}^n$;
but it's now more
natural to think of
the latter---and we will do so in what follows---since
our original universe plays no further role in this discussion.
It may also be worth noting that, unlike in the proof of
Lemma \ref{muflemma}, the second moment method is not strong
enough to give the exponential bound in Claim A.

\mn
{\em Claim B.}
If $X\sub Q$,
$\mu(X) =\gb> \exp[-o(n/\log^2n)]$ and $\gz =o(2^{-k})$,
then for a random $k$-clause $C$,
$$
\Pr(\mu (C\wedge X) >(1-\gz)\mu(X)) < 1/2
$$
(where $C\wedge X =\{x\in X: C\sim x\}$).

\mn
{\em Remark.}  This is probably true for
$\gb$ greater than something like $\exp[-n/k]$.
The bound in the claim
is just what the proof gives, and is more than enough
for us since we're really interested in much larger
$\gb$.

\medskip
To see that Claim B implies Claim A,
set $f_j=\wedge_{i=1}^jC_i$ and notice that $\mu(f)\geq \gb$
implies (for example)
\beq{imufi}
|\{i: \mu(f_i) <(1-\tfrac{5\ln(1/\gb)}{m})\mu(f_{i-1})\}| < m/5
\enq
(and, of course, $\mu(f_i)\geq \gb ~\forall i$).
But if we take
$\gb =\exp[-n^{\gd}]$ then our choice of parameters gives
$$\gz:=5m^{-1}\ln(1/\gb) =o(2^{-k})$$
(using $m/\ln (1/\gb) =
\Theta(n^{(1+\gd)\log (1+2\gd)+1-\gd})$ and
$2^k = n^{1+\gd}$),
so Claim B bounds the probability of \eqref{imufi} by
$$
\Cc{m}{m/5} 2^{-4m/5} =o(2^{-\nn}).
$$\qed

\mn
{\em Proof of Claim B.}
Let $G$ be the bipartite graph on
$Q\cup W$, where $W$ is the set of $(n-k)$-dimensional
subcubes of $Q$ and, for $(x,D)\in Q\times W$,
we take $x\sim D$ if $x\in D$.
(So we've gone to complements:  for a clause $C$
the corresponding subcube is $D=\{y:C\not\sim y\}$,
so $C\wedge X =X\sm D$.)

Assuming Claim B fails at $X$, fix $\SSS \sub W$ with $|\SSS |=|W|/2$ and
$$D\in \SSS
  ~\Ra ~\mu(D\cap X) < \gz\mu(X),$$
and set $\TTT =W\sm \SSS $.

Consider the experiment:
(i) choose $x$ uniformly from $X$;
(ii)  choose $D$ uniformly from the members of $W$ containing $x$;
(iii)  choose $y$ uniformly from $D$.

\mn
{\em Claim C.}  $\Pr(y\in X)> (2-o(1))\gb$.

\mn
{\em Proof.}
Since each triple $(x,D,y)$ with $x\in X$ and $x,y\in D$ is
produced by (i)-(iii) with probability
$|X|^{-1}\cdot 2^k|W|^{-1}\cdot 2^{k-n}$,
we just need to show that the number of such triples with
$y\in X$ is at least
$$
(2-o(1))\gb |X||W|2^n2^{-2k}  = (2-o(1))|X|^2|W|2^{-2k}.
$$
Writing $d$ for degree in $G$, we have
$$
\sum_{x\in X}d_\SSS(x) =\sum_{D\in \SSS }d_X(D) < |\SSS |\gz
 |X|,
$$
implying
\begin{eqnarray*}
\sum_{D\in \TTT }d_X(D)& = &
\sum_{x\in X}(d(x)-d_\SSS (x))\\
&>& |X||W|2^{-k} - \gz |\SSS ||X| = (1-o(1))|X||W|2^{-k}.
\end{eqnarray*}
The number of $(x,D,y)$'s as above is thus
\begin{eqnarray*}
\sum_{D\in W} d_X^2(D)
&\geq &
\sum_{D\in \TTT } d_X^2(D)
\geq
(\sum_{D\in \TTT } d_X(D))^2/|\TTT |\\
&> &(1-o(1))|X|^2|W|^2|\TTT |^{-1}2^{-2k}
= (2-o(1))|X|^2|W|2^{-2k}.
\end{eqnarray*}
\qed

Let $T(x)$ be
the random element of $Q$ gotten from $x$
by choosing $K$ uniformly from $\C{[n]}{k}$ and
randomly (uniformly, independently) revising the $x_i$'s with $i\not\in K$.
Then $y$ gotten from $x$ by (ii) and (iii) above is just $T(x)$,
so
the next assertion contradicts Claim C, completing the proof of
Claim B (and Theorem \ref{prop2}).

\mn
{\em Claim D.}
If $\mu(X) > \exp[-o(n/\log^2n)]$
and $x$ is uniform from $X$, then $\Pr(T(x)\in X) < (1+o(1))\mu(X)$.

\mn
{\em Remark.}
If $X$ is a subcube of codimension $n/k$, say
$X=\{x:x \equiv 0 ~\mbox{on}~L\}$ with $|L |=n/k$,
then for any $x\in X$,
$$
\Pr(T(x)\in X) =\sum_t\Pr(|K\cap L |=t)2^{-(|L|-t)}
=\mu(X)\sum_t\Pr(|K\cap L |=t)2^t,
$$
and, since $|K\cap L|$ is essentially Poisson with mean 1,
the sum is
approximately
$ e^{-1}\sum_t2^t/t! = e$.
So Claim D fails for $\mu(X) = 2^{-n/k}$ and, as earlier,
it's natural to guess that it holds if $\mu(X)$ is much
bigger than this.

\mn
{\em Proof of Claim D.}
Let $Q_r=\{y\in Q:|y|\leq r\}$.
The assumption on $\mu(X)$ implies that
$\mu(Q_{r-1})<\mu(X) \leq \mu(Q_r)$
for some $r >(1/2-o(1/k))n$, so Claim D follows from

\mn
{\em Claim E.}
If $\vp =o(1/k)$ and $r > (1/2-\vp)n$, then for any $x\in Q$ and
$X\sub Q$ with
\beq{muX}
\mu(X) \leq \mu(Q_r),
\enq
\beq{PrTx}
\Pr(T(x)\in X) < (1+o(1))\mu(Q_r).
\enq

\mn
{\em Proof}.
We may assume $x=\underline{0}$, so that $\Pr(T(x)=y)$ is a decreasing
function of $|y|$.
We thus maximize $\Pr(T(x)\in X)$
subject to \eqref{muX} by taking
$
X=Q_r,
$
and \eqref{PrTx} is then a routine calculation
using
$$\mu(X) =\Pr({\rm Bin}(n,1/2)\leq r)$$
and
$$\Pr(T(x)\in X) =\Pr({\rm Bin}(n-k,1/2)\leq r)$$
(where ${\rm Bin}(\cdot,\cdot)$ denotes a binomially
distributed r.v.).\qed

\mn
{\bf Acknowledgment}
We would like to thank Roy Meshulam for helpful conversations.
Part of this work was carried out while the first author was visiting Jerusalem under BSF grant 2006066.

\bn
Department of Mathematics\\
Rutgers University\\
Piscataway NJ 08854\\
jkahn@math.rutgers.edu\\

\bn
Institute of Mathematics\\
Hebrew University of Jerusalem,\\
Jerusalem, Isreal\\
and\\
Department of Mathematics\\
Yale University\\
New Haven, CT\\
kalai@math.huji.ac.il

\end{document}